\documentclass{amsart}
\usepackage{graphicx}
\usepackage{amsmath}
\usepackage{amsfonts}
\usepackage{amssymb}
\newtheorem{theorem}{Theorem} [section]
\newtheorem{thm}[theorem]{Theorem}

\newtheorem{example}[theorem]{Example}

\newtheorem{prop}[theorem]{Proposition}

\begin{document}

\title{Milnor's isotopy invariants and generalized link homotopy}

\author{Thomas Fleming$^{1}$ \and Akira Yasuhara}

\begin{abstract}

It has long been known that a Milnor invariant with no repeated index is an invariant of link homotopy.  We show that Milnor's invariants with repeated indices are invariants not only of isotopy, but also of self $C_{k}$-moves. A self $C_{k}$-move is a natural generalization of link homotopy based on certain degree $k$ clasper surgeries, which provides a filtration of link homotopy classes.

\end{abstract}

\maketitle

\section{Introduction}

In his landmark 1954 paper \cite{m1}, Milnor introduced his epinonomous higher order linking numbers, and proved that when the multi-index $I$ had no repeated index, the $\mu(I)$ were invariants of link homotopy.  In a follow up paper \cite{m2}, Milnor explored some of the properties of invariants with repeated indices.  

Milnor's invariants have been studied extensively since 
that time.  For 
string links, these invariants are known to be finite type \cite{bn, lin}, and in fact related to (the tree part of) the Kontsevich integral in a natural and beautiful way \cite{habmas}. \footnotetext[1]{The first author was supported by a Japan Society for the 
Promotion of Science Post-Doctoral Fellowship for Foreign Researchers, 
(Short-Term), Grant number PE05003}

By work of Habiro \cite{ha}, the finite type invariants of knots are intimately related to claper surgery.  Taking the view that link homotopy is generated by degree one clasper surgery on the link where both leaves of the clasper are on the same component, is it possible that Milnor's isotopy invariants have some relation with clasper surgery?

The answer is yes.  Let us define a \emph{self $C_{k}$-move} on a link $L$ to be a degree $k$ simple tree clasper surgery on $L$ where all leaves of the clasper are on the same component.  If $L'$ is obtained from $L$ by a sequence of self $C_{k}$-moves, we call $L$ and $L'$ \emph{self $C_{k}$-equivalent}.  These moves were introduced in \cite{shiyas}, and as we have mentioned, self $C_{1}$-equivalence is link homotopy.  

For an $n$-component link, Milnor invariants are specified by a 
multi-index $I$, where the entries of $I$ are chosen from $\{1, \ldots 
n\}$.  
Let $r(I)$ denote the maximum number of times that any index appears. A
Milnor invariant $\mu(I)$ is called \emph{realizable} if there exists a 
link $L$ with $\mu_{L}(I) \neq 0$.

Our main result is the following.

\medskip

\textbf{Theorem \ref{mainthm}} \emph{Let $\mu(I)$ be a realizable Milnor 
number.  Then $\mu(I)$ 
is an invariant of self $C_{k}$-equivalence if and only if $r(I) \leq k$.}
\medskip

Notice that a self $C_{k}$-equivalence can be realized by self ${C_{k'}}$-moves when $k' < k$, and thus self $C_{k}$-equivalence classes form a filtration of link homotopy classes.  Moving to larger and larger $k$ provides more and more information about the structure of isotopy classes of links.  

The classification of links up to link homotopy has been completed by Habegger and Lin \cite{hl}.  However, very little is known about the structure of links under these higher order moves.  Nakanishi and Ohyama have classified two component links up to self $C_{2}$-equivalence \cite{no1, no2, no3}, and it is known that boundary links are self $C_{2}$-equivalent to the trivial link \cite{yasuhara}.  

It is well known that Milnor's link homotopy invariants vanish if and only if the link is link homotopic to the unlink.  However, in Example \ref{c3eg} we will produce a boundary link that is not self $C_{3}$-equivalent to the unlink.  As all Milnor invariants vanish for boundary links, this example demonstrates that self $C_{k}$-equivalence behaves very differently for $k>1$.  Much work remains to be done before we have a clear understanding of self $C_{k}$-equivalence.

\section{Invariance under self $C_{k}$-moves}

We define a self $C_{k}$-move to be surgery on a simple, tree-like degree $k$ clasper, all of whose leaves are on the same component.  Similarly, $L$ is self $C_{k}$-equivalent to $L'$ if $L'$ can be obtained from $L$ by a sequence of self $C_{k}$-moves.    

Note that link homotopy is the same as self $C_{1}$-equivalence, and the self-delta equivalence of Shibuya \cite{shi} is the same as self $C_{2}$-equivalence.  

For an $n$-component link, Milnor invariants are specified by a 
multi-index $I$, where the entries of $I$ are chosen from $\{1, \ldots 
n\}$.  
Let $r(I)$ denote the maximum number of times that any index appears.  A 
Milnor invariant $\mu(I)$ is called \emph{realizable} if there exists a 
link $L$ with $\mu_{L}(I) \neq 0$.

\begin{thm} Let $\mu(I)$ be a realizable Milnor number.  Then $\mu(I)$ is 
an invariant of self $C_{k}$-equivalence if and only if $r(I) \leq k$.
\label{mainthm}
\end{thm}

\begin{proof}

Suppose that $r(I) \leq k$ and $L$ is self $C_{k}$-equivalent to $L'$. Since Milnor's invariants are well known to be isotopy invariants, it suffices to check a single self $C_{k}$-move realized by a clapser $c$.

The case $k=1$ is that of link homotopy, and was proven by Milnor in \cite{m1}.

Thus, we need only consider the case $k \geq 2$. 
To calculate $\mu_{L}(I)$ (resp. $\mu_{L'}(I)$) where $I$ has repeated indices, we may instead study the link homotopy invariants of $\overline{L}$ ($\overline{L'}$), the link obtained by taking the appropriate number of zero framed parallels of the components of $L$ ($L'$) \cite{m2}. In particular, if component $i$ of $\overline{L}$ corresponds to componet $h(i)$ of $L$, then $\mu_{\overline{L}}(i_{1},i_{2} \ldots i_{m}) = \mu_{L}(h(i_{1}),h(i_{2}) \ldots h(i_{m}))$.

Since $k \geq 2$, the self $C_{k}$-move preserves the framing of the components of $L$.  That is, the clasper surgery carrying $L$ to $L'$ carries $\overline{L}$ to $\overline{L'}$.

The self $C_{k}$-move is realized by surgery on $c$, that is, a clasper with $k+1$ leaves where all of these leaves land on the same component of $L$.  Replace $L$ with $\overline{L}$. 

As $r(I) \leq k$, $\overline{L}$ has at most $k$ copies each component of $L$, so each leaf of the clasper $c$ grips at most $k$ components of $\overline{L}$.  Further, as each leaf of $c$ was on the same component of $L$, each leaf of $c$ holds the same components of $\overline{L}$.

Use the zip construction on $c$ to produce simple claspers on 
$\overline{L}$. For a version of the zip construction without ``boxes'', 
see 4.2 of \cite{coteich}.  

Each simple clasper has at least two leaves on 
the same 
component of $\overline{L}$, and so is realized by link homotopy. See Proposition \ref{easyprop}.

We have shown that $\overline{L}$ is link homotopic to $\overline{L'}$, so  $\mu_{\overline{L}} = \mu_{\overline{L'}}$ for all link homotopy invariants. Hence, when $r(I) \leq k$,  $\mu_{L}(I) = \mu_{L'}(I)$.  

\medskip

Suppose that $r(I) > k$.  Since the Milnor invariant is 
realizable, by work of Cochran \cite{cochran}, there exists a link $L$ 
where $\mu_{L}(I) \neq 0$, and all lower order Milnor invariants vanish. Such a link can be obtained by Bing doubling a Hopf link along a trivalent tree whose leaves are labeled by $I$, and then band summing components with the same label.  This process is the same as a simple tree clasper surgery on the unlink, using the same tree, where leaves labeled $i$ grasp the $i$th component.  

Since some index $i$ in $I$ is repeated at least $k+1$ times, this clasper has at least $k+1$ leaves on component $i$ and so can be realized by self $C_{k}$-moves. See Proposition \ref{easyprop}.

As all Milnor invariants of the unlink vanish, we have that $\mu_{L}(I)$ is not a self $C_{k}$ invariant when $r(I) > k$.

\end{proof}

We end this section with an example that demonstrates the power of our methods.

Self $C_{2}$-equivalence of links is equivalent to self-delta equivalence, a relation that has been extensively studied.   Much of this work relies on invariants of self $C_{2}$-(self-delta)-equivalence that are based on the coefficients of the Conway polynomial.  

\begin{figure}[hbtp]
\centering
\begin{picture}(140,122)
\includegraphics{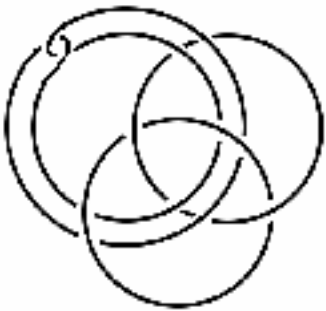}
\end{picture}
\caption{The Whitehead double of the Borromean rings.}
\label{whbring}
\end{figure}

Let $L$ be the Whitehead double of the Borromean rings (equivalently the 
Bing double of the Whitehead link). See Figure \ref{whbring}.  The 
Alexander polynomial of this link is trivial\footnote{J. Hillman has 
pointed out that the three variable Alexander polynomial of $L$ is also 
$0$.}, but $\mu_{L}(123123)=1$.  
Thus, $L$ has trivial Conway polynomial, but is not self 
$C_{2}$-equivalent to the unlink.  

\section{Further results}

In the proof of Theorem \ref{mainthm} we used the fact that simple clasper surgery with multiple leaves on a single component could be converted to a self $C_{k}$-equivalence.  In this section, we will find a partial converse to that fact and demonstrate its usefulness in the study of self $C_{k}$-equivalence.

We say $L$ and $L'$ are $C_{m}^{k}$-equivalent if $L$ can be transformed into $L'$ by ambient isotopy and degree $m$ simple tree clasper surgery, where at least $k$ leaves of the clasper grasp component $i$.

We may now state the fact mentioned above precisely.

\begin{prop} If $L$ is $C_{m}^{k+1}$-equivalent to $L'$ then $L$ is self $C_{k}$-equivalent to $L'$.
\label{easyprop}
\end{prop}

\begin{figure}[hbtp]
\centering
\begin{picture}(380,176)
\includegraphics{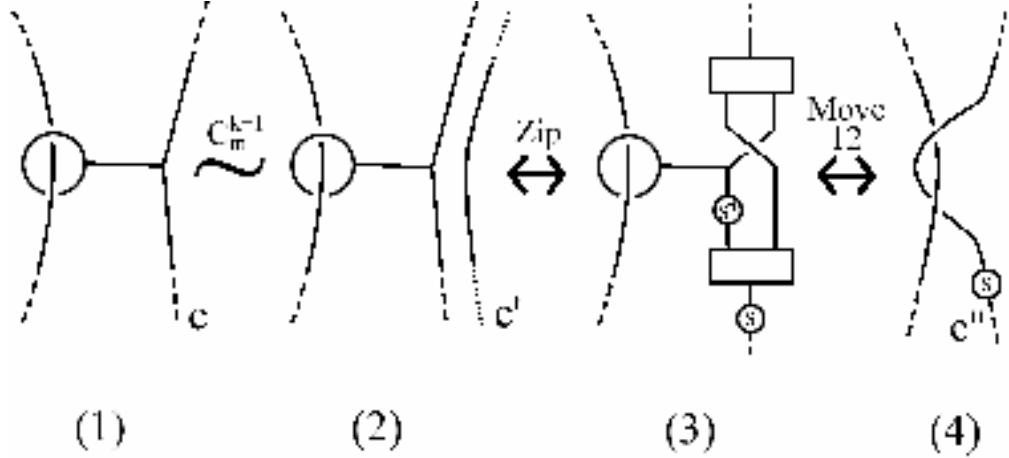}
\end{picture}
\caption{Reducing a $C_{m+1}^{k+1}$ clasper to $C_{m}^{k+1}$ claspers. The 
symbol $s$ ($s^{-1}$) represents a positive (negative) half twist.}
\label{clasperslide}
\end{figure}

\begin{proof}

We will work by induction.  The base case is $m=k$, and is trivial. 

Suppose now that $c$ is a clasper representing a $C_{m+1}^{k+1}$ move.  We will show that it is $C_{m}^{k+1}$-equivalent to a $C_{m}^{k+1}$ move.  Begin with the clasper $c$.  Introduce a clasper $c'$ representing a $C_{m}^{k+1}$-move as in Figure \ref{clasperslide} (2).  We choose $c'$ so that by reversing the zip construction, we obtain the clapser with boxes shown in Figure \ref{clasperslide}.  By Move 12 of \cite{ha}, we obtain the clasper $c''$, which is degree $m$ and has $k+1$ feet on component $i$.  Thus, we can reduce any $C_{m}^{k+1}$-move to a sequence of $C_{k}^{k+1}$-moves, which is merely a self $C_{k}$-equivalence.


Figure \ref{clasperslide} demonstrates the necessary moves. One may find it easier to read this figure from right to left.  
\end{proof}

Proposition \ref{easyprop} is a useful, though elementary, application of the properties of claspers.  Less obvious, however, is that in certain circumstances we have a converse.  

\begin{prop}Let $L$ be an $n$ component Brunnian link, and $U$ the $n$ component unlink.  Then $L$ is self $C_{k}$-equivalent to $U$ if and only if $L$ is $C_{n+k-1}^{k+1}$-equivalent to $U$.
\label{bruneq}
\end{prop}

\begin{figure}[hbtp]
\centering
\begin{picture}(382,146)
\includegraphics{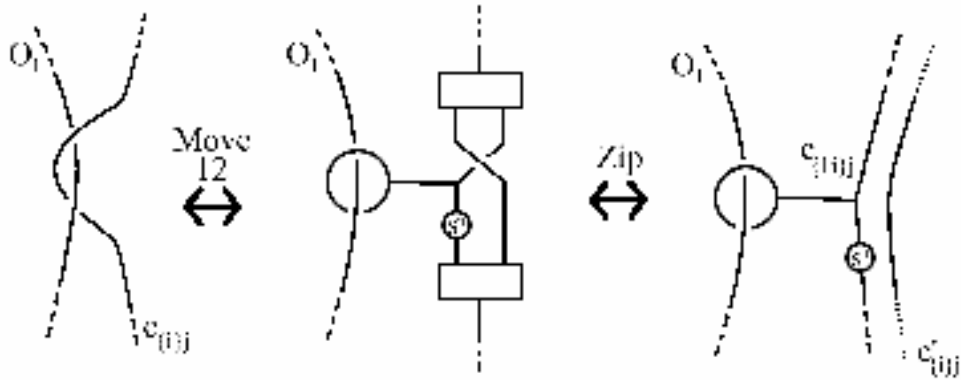}
\end{picture}
\caption{Altering the clasper $c_{\{i\}j}$ to reflect a crossing change between that clasper and $O_{1}$.}
\label{clasper}
\end{figure}

\begin{proof}

The `if' part follows from Proposition \ref{easyprop}.  We will show the `only if' part.  The idea of this proof is similar to that of Theorem 1.2 in \cite{miya}.  While in \cite{miya} they use the `\emph{band description}' defined in \cite{taniyas}, we will use clasper theory.

Assume that $L$ and $U$ are self $C_k$-equivalent. 
Let $L_i$ be the $i$th component of $L$. 
We can describe $L$ as the unlink 
$O_1\cup\cdots \cup O_n$ with claspers $c_{\{i\},j}$ representing 
the self $C_k$-moves on  the $i$th component. 
Since $L$ is Brunnian, $L\setminus L_1$ is the unlink. 
Thus, $L$ is obtained from the trivial link 
by $C_k$-moves corresponding to $c_{\{1\},j}$ and 
by crossing changes 
between $O_1$ and $c_{\{i\},j}~(i\geq 2)$. 

Passing a clasper through $O_{1}$ results in two new claspers, one of 
degree $k+1$ which has a new leaf on $O_{1}$, and one of degree $k$ which does not.  Call the former $c_{\{1i\},j}$, and the later $c_{\{i\},j}'$.  See Figure \ref{clasper}.

We can now express $L$ as the unlink $O_1\cup\cdots \cup O_n$ with claspers $c_{\{1\},j}$ of degree $k$, claspers $c_{\{i\},j}'~(i\geq 2)$ of degree $k$ and 
claspers $c_{\{1i\},j}~(i\geq 2)$ of degree $k+1$.  

Let $L'$ be the link obtained by surgering the unlink $O_1\cup\cdots \cup O_n$ only on the $c_{\{i\},j}'$. By the construction of the $c_{\{i\},j}'$, $L' = L \setminus L_{1} \coprod O_{1}$.  Since $L$ is Brunnian, $L'$ is the unlink, and so the claspers $c_{\{i\},j}'$ act trivially.    

Thus we can express $L$ as the unlink  
$O_1\cup\cdots \cup O_n$ with claspers $c_{\{1\},j}$ of degree $k$ and 
claspers $c_{\{1i\},j}~(i\geq 2)$ of degree $k+1$ that have 
leaves on $O_1\cup O_i$.

Since $L\setminus L_2$ is the unlink, repeating the argument above 
shows that $L$ is expressed as the unlink  
$O_1\cup\cdots \cup O_n$ with claspers $c_{\{12\},j}$ of degree $k+1$ 
that have leaves on $O_1\cup O_2$
and 
claspers $c_{\{12i\},j}~(i\geq 3)$ of degree $k+2$ that have 
leaves on $O_1\cup O_2\cup O_i$. 

Repeating these process for each $L_i$, finally, we can express 
$L$ as the unlink  
$O_1\cup\cdots \cup O_n$ with claspers $c_{\{12,...,n\},j}$ of degree $k+n-1$ 
that have leaves on $O_1\cup O_2\cup\cdots\cup O_n$. Hence $L$ is 
$C_{k+n-1}$-equivalent to the unlink.  
\end{proof}

\begin{figure}[hbtp]
\centering
\begin{picture}(377,133)
\includegraphics{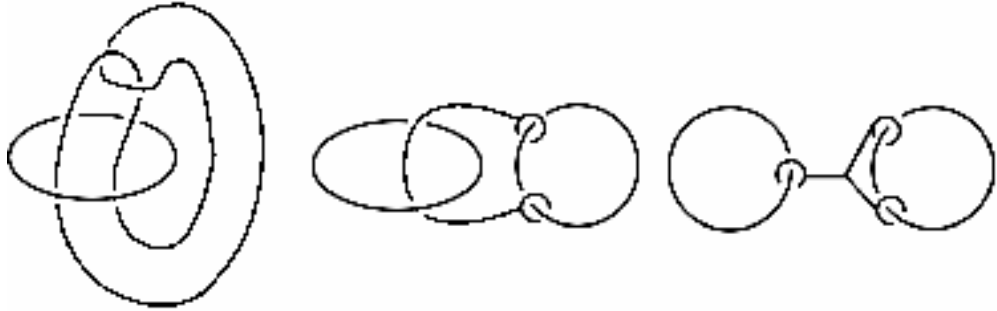}
\end{picture}
\caption{Proposition \ref{bruneq} for the Whitehead link.}
\label{whitehead}
\end{figure}

Figure \ref{whitehead} illustrates the argument of Proposition 
\ref{bruneq} for the Whitehead link.  The Whitehead link is Brunnian and 
link homotopic to the unlink. The center image of Figure \ref{whitehead} shows it as the unlink with one self $C_{1}$ clasper.  The right hand image shows the $C_{2}$-equivalence of the Whitehead link and the unlink.  

\begin{figure}[hbtp]
\centering
\begin{picture}(159,164)
\includegraphics{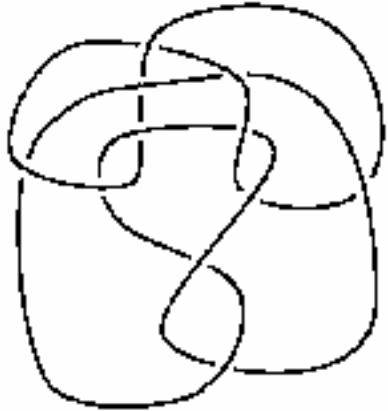}
\end{picture}
\caption{The Hopf link with both components Whitehead doubled.}
\label{whwhhopf}
\end{figure}

\begin{example}

The link in Figure \ref{whwhhopf} is obtained from Whitehead doubling both 
components of the Hopf link and is a boundary link, so all Milnor 
invariants vanish.  However, examine the Jones polynomial.  $$J(L) = 
q^{-\frac{9}{2}} -2q^{-\frac{7}{2}} +q^{-\frac{5}{2}}-q^{-\frac{3}{2}} 
-q^{\frac{3}{2}} + q^{\frac{5}{2}} -2 q^{\frac{7}{2}} + q^{\frac{9}{2}}$$ 
We obtained the Jones polynomial from the Knot Atlas 
\footnote{http://www.math.toronto.edu/$\sim$drorbn/KAtlas, link L11n247} 
and by calculation using Knot 
\footnote{http://www.math.kobe-u.ac.jp/$\sim$kodama/knot.html}. Evaluating 
the third derivative of $J(L)$ at 1, we obtain a value distinct from that 
of the third derivative of Jones polynomial of the unlink at 2.  Since 
this is a finite type invariant of degree three, $L$ is not $C_{4}$-equivalent to the trivial link.  Using Proposition \ref{bruneq} we see 
that it is not self $C_{3}$-equivalent to the trivial link.  In fact, 
since $L$ is Brunnian, if $L$ is self $C_{k}$-equivalent to a split link, 
it is self $C_{k}$-equivalent to the trivial link. 
Thus, $L$ is not self $C_{3}$-equivalent to a split link.  \label{c3eg} 
\end{example}

Example \ref{c3eg} is interesting for other reasons.  It is well known 
that Milnor's link homotopy invariants vanish if and only if $L$ is link 
homotopic to the unlink.  The example above shows that a similar statement 
is not true for Milnor's self $C_{3}$ invariants.  While the link in 
example is a boundary link, it is not self $C_{3}$-equivalent to the 
trivial link 
or even a split link.  It would be interesting to know what the vanishing 
of Milnor's self $C_{k}$ invariants implies about the self 
$C_{k}$-equivalence class of $L$.

For two component links, vanishing $\mu_{L}(12)$ and $\mu_{L}(1122)$ implies that the link is self $C_{2}$-equivalent to the unlink \cite{no3}.  Shibuya and the second author have recently shown that all boundary links are self $C_{2}$-equivalent to the unlink \cite{yasuhara}, but whether the vanishing of Milnor's self $C_{2}$ invariants is sufficient to show that a link is self $C_{2}$-equivalent to the unlink remains an open question.  


\textsc{University of California San Diego, Department of Mathematics, 9500 Gilman Dr., La Jolla,
 CA 92093-0112, United States}

\emph{E-mail address:} \texttt{tfleming@math.ucsd.edu}

\medskip

\textsc{Tokyo Gakugei University, Department of Mathematics, Koganei-shi, Tokyo 184-8501, Japan}

\emph{E-mail address:} \texttt{yasuhara@u-gakugei.ac.jp}

\end{document}